%% file: panda.tex
\documentclass[letterpaper, 10 pt, conference]{ieeeconf}  

\IEEEoverridecommandlockouts                              
\usepackage{graphics} 
\usepackage{epsfig} 
\usepackage{times} 

\usepackage{amsmath}
\usepackage{color}
\usepackage{pgfplots} 
\usepackage{pgf}
\usepackage{tikz}
\usetikzlibrary{arrows,automata,patterns,chains}
\usepackage{amssymb}
\usepackage{mathtools}
\usepackage{algorithm}
\usepackage{textcomp}
\usepackage{algpseudocode}
\usepackage{amsfonts}
\usepackage{subfiles}
\usepackage{float}
\usepackage{tkz-graph}

\makeatletter
\renewcommand*\env@matrix[1][*\c@MaxMatrixCols c]{%
  \hskip -\arraycolsep
  \let\@ifnextchar\new@ifnextchar
  \array{#1}}
\makeatother  

\newtheorem{assumption}{Assumption}
\newtheorem{theorem}{Theorem}

\newfloat{algorithm}{t}{lop}

\makeatletter
\newcommand\fs@spaceruled{\def\@fs@cfont{\bfseries}\let\@fs@capt\floatc@ruled
  \def\@fs@pre{\vspace{5pt}\hrule height.8pt depth0pt \kern2pt}%
  \def\@fs@post{\kern2pt\hrule\relax}%
  \def\@fs@mid{\kern2pt\hrule\kern2pt}%
  \let\@fs@iftopcapt\iftrue}
\makeatother

\title{\LARGE \bf
PANDA: A Dual Linearly Converging Method for Distributed Optimization over Time-Varying Undirected Graphs}

\author{Marie Maros and Joakim Jald\'{e}n
\thanks{This project has received funding from the European Research Council (ERC) under the European Union's Horizon 2020 research and innovation programme (grant agreement No 742648)}
\thanks{The authors are with the Department of Information Science and Engineering, EECS, KTH, Stockholm, Sweden.
        {\tt\small mmaros@kth.se, jalden@kth.se}}%
}

\begin{document}

\maketitle
\thispagestyle{empty}
\pagestyle{empty}

\begin{abstract}
In this paper we consider a distributed convex optimization problem over time-varying networks. We propose a dual method that converges R-linearly to the optimal point given that the agents' objective functions are strongly convex and have Lipschitz continuous gradients. The proposed method requires half the amount of variable exchanges per iterate than methods based on DIGing, and yields improved practical performance as empirically demonstrated.
\end{abstract}

\section{Introduction}

Solving optimization problems in a distributed manner has become relevant in many engineering applications. These include power system control \cite{giannakis}, distributed resource allocation \cite{tutorial}, distributed estimation \cite{application,estimation}, statistical inference and learning \cite{learning} among others \cite{localization,formation}. In these set-ups partial information of the problem is available to different agents in the network. Further, having a fusion center in which all information is gathered may not be affordable for several reasons. In some instances, the nodes may have gathered massive amounts of data, making the communication of it expensive. In others, the nodes may have to communicate via a wireless network resulting in communication costs. Additionally, due to the random nature of the wireless channel some communication links may fail. Therefore,  procedures that allow for the agents to collaboratively solve optimization problems while exchanging limited amounts of information via an unreliable network are required.

In this paper we focus on the distributed optimization problem 
\begin{align}
\label{eq:original}
& \underset{\mathbf{\bar{x}} \in \mathbb{R}^p}{\text{min}} \quad \bar{f}(\mathbf{\bar{x}}) = \sum_{i=1}^n f_i(\mathbf{\bar{x}}),
\end{align}
where each function $f_i$ is known exclusively to a single agent of the network. The agents' goal is to collaboratively by communicating only with their direct neighbors find a solution to \eqref{eq:original} without explicitly exchanging the functions $f_i.$ We will represent the network connecting the different agents via a time varying graph. This allows for applications such as sensor and vehicle networks in which some communication links may come and go.

Many methods that solve problems of the form of \eqref{eq:original} have been proposed in recent years. First order primal methods typically require diminishing step-sizes to converge to a solution even if the functions $f_i$ are strongly convex and have Lipschitz continuous gradients \cite{nedic}, \cite{decentralized_gradient}. While both \cite{nedic} and \cite{decentralized_gradient} have been shown to converge given time-varying networks, selecting a diminishing step-size will lead to slow convergence rates. Dual and primal-dual methods typically converge to the optimal solution using a constant step-size \cite{extra,pextra,harnesing,optimal,dual_decomposition,admm}. However, extensions of these methods that converge linearly under time varying graphs have not been proposed. This said, several methods that converge under time varying graphs have been recently proposed \cite{random_networks,interpret,next,non_convex}. In \cite{random_networks} the authors propose a primal-dual algorithm that achieves sub-linear convergence rates. By proposing a different error correction scheme than that in \cite{extra}, the authors of \cite{interpret} propose the decentralized inexact gradient tracking (DIGing), a method which converges linearly to the optimal solution even when the graph is time varying. In-network nonconvex optimization (NEXT) \cite{next} is a method to solve non-convex problems of the form \eqref{eq:original} on undirected graphs.  In \cite{non_convex} the authors propose the successive convex approximation over time-varying digraphs (SONATA) to solve non-convex optimization problems over time-varying directed graphs. SONATA can be particularized to yield the same iterates as Push-DIGing \cite{interpret} which modifies DIGing to be applicable to directed graphs.
Under the assumption that the functions $f_i$ are strongly convex and have Lipschitz continuous gradients NEXT, SONATA and DIGing will converge linearly to the optimal point. Further, the three of them are primal methods that rely on tracking the average primal variable and the average gradient. This procedure, while effective, requires that the agents exchange the primal variables and the gradients at each iterate. 

In this paper, we propose \underline{p}rimarily \underline{a}veraged \underline{n}etwork \underline{d}ual \underline{a}scent (PANDA), a novel linearly converging dual ascent-based method for distributed optimization over time varying undirected graphs. PANDA can be shown to converge linearly to the optimal point of \eqref{eq:original} given that the $f_i$ are strongly convex and have Lipschitz continuous gradients. When compared to DIGing, PANDA requires only the exchange of the primal variables, at the expense of iterates that are computationally more expensive. In terms of rate, PANDA performs better in practice as shown in our experiments. Further, in \cite{optimal} the authors establish that, given a static connected network, a Nesterov accelerated version of dual ascent achieves optimal convergence rates. This further motivates the proposal of a dual method that can deal with time-varying graphs. The theoretical analysis of the convergence properties of an accelerated version are left for future work. However, we verify experimentally that we obtain a convergence gain by applying Nesterov's acceleration scheme.

The remaining of the paper is structured as follows. In Section \ref{section:algorithm} we provide an intuition for the derivation of PANDA and give its formal introduction. After that, in Section \ref{section:convergence}, we formally present the assumptions under which this paper's main theorem holds. Then, the main theorem containing the algorithm's convergence result is stated. In Section \ref{section:sketch} we provide a rough sketch of the proof of the main theorem. Details are omitted due to space constraints. In order to experimentally verify the findings of this paper we provide empirical evidence of the convergence of PANDA in Section \ref{section:numerical}. Finally, the paper ends with some concluding remarks in Section \ref{section:conclusions}.

\section{Algorithm and algorithm intuition \label{section:algorithm}}
In this section we will lay out the intuition for PANDA. We start by reformulating \eqref{eq:original} in the usual way corresponding to solving it using distributed dual ascent \cite{optimal, dual_decomposition,tutorial}. Once we have done this, we will observe the main inconveniences with dealing with the standard reformulation when we are solving the problem over a time-varying graph. Avoiding these inconveniences will lead us to a different problem formulation that relies on the fully connected network. Finally, by approximating the effect of a fully connected network we obtain the desired algorithm.

To this end, let the nodes $i =1,\hdots,n$ be connected via the undirected graph $\mathcal{G}(\mathcal{V},\mathcal{E})$ with vertices $\mathcal{V} \triangleq \{1,\hdots,n\}$ and edges $\mathcal{E}.$ Two nodes $i$ and $j$ can only communicate directly with each other if the edge $(i,j) \in \mathcal{E}.$ Further, let $\mathbf{U}$ denote the graph's communication matrix with the following properties.
\begin{assumption}[Communication Matrix $\mathbf{U}$ \cite{optimal}] \label{assumption:communication} The matrix $\mathbf{U}$ fulfills the following properties:
\begin{enumerate}
\item [(P1)] $\mathbf{U}$ is an $n \times n$ symmetric matrix,
\item [(P2)] $\mathbf{U}$ is positive semi-definite,
\item [(P3)]$\text{null}\{\mathbf{U}\} = \text{span}(\mathbf{1}_n),$ where $(1,\hdots,1)^T = \mathbf{1}_n \in \mathbb{R}^n,$
\item [(P4)]$\mathbf{U}$ is defined on the edges of the network, i.e., $u_{ij} \neq 0$ only if $(i,j) \in \mathcal{E}.$
\end{enumerate}
\end{assumption}

Under Assumption \ref{assumption:communication} the problem in \eqref{eq:original} can be equivalently written as
\begin{subequations}
\label{eq:standard_dual_reformulation}
\begin{align}
& \underset{\mathbf{x} \in \mathbb{R}^{np}}{\text{min}} \quad f(\mathbf{x}) = \sum_{i=1}^n f_i(\mathbf{x}_i) \label{eq:standard_dual_reformulation:objective} \\
& \text{s.t.} \quad \left(\mathbf{U}^{1/2} \otimes \mathbf{I}_p \right)\mathbf{x} = \mathbf{0}, \label{eq:standard_dual_reformulation:constraint}
\end{align}
\end{subequations}
where $\mathbf{x}^T = (\mathbf{x}_1^T,\hdots,\mathbf{x}_n^T),$ $\mathbf{I}_p \in \mathbb{R}^{p \times p}$ denotes the identity matrix and $\mathbf{U}^{1/2}$ is a symmetric matrix such that $\mathbf{U}^{1/2}\mathbf{U}^{1/2} = \mathbf{U}.$
Note that each vector $\mathbf{x}_i, \, i = 1,\hdots,n$ denotes the local copy of the variable $\mathbf{\bar{x}}$ (c.f. \eqref{eq:original}) held by node $i.$ Together, all local copies $\mathbf{x}_i$ constitute the global variable $\mathbf{x} \triangleq (\mathbf{x}_1^T,\hdots,\mathbf{x}_n^T)^T.$ The constraint \eqref{eq:standard_dual_reformulation:constraint} enforces $\mathbf{x}$ to be \emph{consensual}, i.e., enforces that all nodes agree $\mathbf{x}_1 = \mathbf{x}_2 =\hdots = \mathbf{x}_n.$  Problems \eqref{eq:original} and \eqref{eq:standard_dual_reformulation} are equivalent in the sense that if $\bar{\mathbf{x}}^{\star}$ solves \eqref{eq:original}, then $\mathbf{x}^{\star} = (\bar{\mathbf{x}}^{\star T},\hdots,\bar{\mathbf{x}}^{\star T})^T$ solves \eqref{eq:standard_dual_reformulation}.

The problem in \eqref{eq:standard_dual_reformulation} can be solved in a distributed fashion using distributed dual ascent by performing the iterates
\begin{subequations}
\label{eq:standard_dual_ascent}
\begin{align}
\label{eq:standard_dual_ascent:obtain}
 \mathbf{x}_i(k+1) := & \quad \underset{\mathbf{x}_i \in \mathbb{R}^p}{\text{min}} \quad f_i(\mathbf{x}_i) - \mathbf{y}_i(k)^T\mathbf{x}_i, \, i=1,\hdots,n \\
\label{eq:standard_dual_ascent:ascent}
 \mathbf{y}_i(k+1) := & \quad \mathbf{y}_i - c \sum_{j \in \mathcal{N}_i \cup \{i\}} u_{ij} \mathbf{x}_j(k+1)
\end{align}
\end{subequations}
where $\mathcal{N}_i \triangleq \{j:(j,i) \in \mathcal{E}\}$ denotes the neighborhood of $i,$ $c > 0$ is an appropriately selected step-size and $u_{ij}$ is the weight node $i$ assigns to the information coming from node $j.$ The quantity $u_{ij}$ corresponds to the element $(i,j)$ of matrix $\mathbf{U}.$ Note that in order to perform the dual ascent step \eqref{eq:standard_dual_ascent:ascent} the nodes have to exchange their update $\mathbf{x}_i(k+1)$ with their immediate neighbors.

If we now assume that the graph $\mathcal{G}(k) = (\mathcal{V},\mathcal{E}(k))$ is allowed to change over time and to become disconnected in some instances, (P3) in Assumption \ref{assumption:communication} would not be fulfilled anymore. This would break the equivalence between problems \eqref{eq:original} and \eqref{eq:standard_dual_reformulation}. Hence, \eqref{eq:standard_dual_ascent} would be  performing an iterate to solve independent optimization problems across the network
instead of attempting to approximate the behavior of solving \eqref{eq:standard_dual_reformulation}. This problem is essentially due to two factors which we summarize as:
\begin{itemize}
\item the problem formulation in \eqref{eq:standard_dual_reformulation} is graph dependent,
\item there is no averaging effect in \eqref{eq:standard_dual_ascent} across time to take advantage of the structure of the different graphs $\mathcal{G}(k)$.
\end{itemize}

Addressing the first concern can be done by seeing what all connected graphs have in common. For this we write the dual of \eqref{eq:standard_dual_reformulation}
\begin{equation}
\label{eq:dual_connected}
\underset{\boldsymbol{\lambda} \in \mathbb{R}^{np}}{\text{min}} \quad f^*\left(\left(\mathbf{U}^{1/2}\otimes \mathbf{I}_p\right) \boldsymbol{\lambda}\right),
\end{equation}
where 
\begin{equation}
f^*(\mathbf{y}) \triangleq \underset{\mathbf{y} \in \mathbb{R}^{np}}{\text{sup}} \quad \mathbf{y}^T\mathbf{x} - f(\mathbf{x})
\end{equation}
denotes the convex conjugate of $f$ and $\boldsymbol{\lambda}$ denotes the dual multiplier associated to \eqref{eq:standard_dual_reformulation:constraint}. The dual can be re-formulated by introducing an additional variable as 
\begin{equation}
\underset{\mathbf{y} \in \mathbb{R}^{np},\boldsymbol{\lambda} \in \mathbb{R}^{np}}{\text{min}} \quad f^*(\mathbf{y}) \quad \text{s.t.} \,\, \left(\mathbf{U}^{1/2} \otimes \mathbf{I}_p\right)\boldsymbol{\lambda} = \mathbf{y}.
\end{equation}
Since $\mathbf{U}^{1/2}$ spans the orthogonal complement of the consensual vectors, for any underlying connected network the constraint $(\mathbf{U}^{1/2}\otimes \mathbf{I}_p)\boldsymbol{\lambda} = \mathbf{y}$ can be replaced by the constraint $(\boldsymbol{\Pi}_{\mathbf{1}_n}\otimes \mathbf{I}_p )\mathbf{y} = \mathbf{0},$ where $\boldsymbol{\Pi}_{\mathbf{1}_n} \triangleq \frac{1}{n}\mathbf{1}_n\mathbf{1}_n^T,$ 
yielding the optimization problem
\begin{equation}
\label{eq:dual_fully_connected}
\underset{\mathbf{y} \in \mathbb{R}^{np}}{\text{min}} \quad f^*(\mathbf{y}) \quad \text{s.t.} \,\, \left(\boldsymbol{\Pi}_{\mathbf{1}_n}\otimes \mathbf{I}_p\right)\mathbf{y} = \mathbf{0},
\end{equation}
which can be solved using projected gradient descent by performing the iterate
\begin{equation}
\mathbf{y}(k+1) := \left(\boldsymbol{\Pi}_{\mathbf{1}_n}^{\perp}\otimes \mathbf{I}_p\right)\left(\mathbf{y}(k) - c \nabla f^*(\mathbf{y}(k)\right),
\end{equation}
where $\boldsymbol{\Pi}_{\mathbf{1}_n}^{\perp} \triangleq \mathbf{I}_n - \boldsymbol{\Pi}_{\mathbf{1}_n}.$ If $\mathbf{y}(0)$ fulfills $(\boldsymbol{\Pi}_{\mathbf{1}_n}^{\perp} \otimes \mathbf{I}_p)\mathbf{y}(0) = \mathbf{y}(0),$ the iterate can be written as
\begin{equation}
\label{eq:equivalent_full}
\mathbf{y}(k+1) := \mathbf{y}(k) - c \left(\boldsymbol{\Pi}_{\mathbf{1}_n}^{\perp} \otimes \mathbf{I}_p\right) \nabla f^*(\mathbf{y}(k)),
\end{equation}
which corresponds to applying dual ascent to the problem 
\begin{subequations}
\label{eq:complete_dual_reformulation}
\begin{align}
 \underset{\mathbf{x} \in \mathbb{R}^{np}}{\text{min}} & \quad f(\mathbf{x}) \label{eq:complete_dual_reformulation:objective} \\
 \text{s.t.}  \quad & \left(\boldsymbol{\Pi}_{\mathbf{1}_n}^{\perp}\otimes \mathbf{I}_p \right)\mathbf{x} = \mathbf{0}. \label{eq:complete_dual_reformulation:constraint}
\end{align}
\end{subequations}
Note that \eqref{eq:complete_dual_reformulation} implicitly assumes that the underlying graph that connects the nodes is fully connected. This is implied by constraint \eqref{eq:complete_dual_reformulation:constraint}. However, if the graph is connected (but not fully connected) the reformulation \eqref{eq:complete_dual_reformulation}, while not practically useful, fulfills (P1)-(P3) in Assumption \ref{assumption:communication}. 
Hence, the reformulation in \eqref{eq:complete_dual_reformulation} can be interpreted as being ``graph independent" in the sense that from the dual problem formulated in \eqref{eq:dual_connected} for any connected graph  we can go to the equivalent dual problem in \eqref{eq:dual_fully_connected}.

By applying dual gradient ascent to \eqref{eq:complete_dual_reformulation} we obtain the iterates
\begin{subequations}
\label{eq:complete_dual_ascent_distributed}
\begin{align}
& \mathbf{x}_i(k+1) := \text{arg } \underset{\mathbf{x}_i \in \mathbb{R}^{p}}{\text{min}} \quad f_i(\mathbf{x}_i) - \mathbf{y}_i(k)^T\mathbf{x}_i,  \quad i=1,\hdots,n \\
& \mathbf{y}_i(k+1) :=  \mathbf{y}_i(k) - c \left( \mathbf{x}_i(k+1) - \frac{1}{n}\sum_{j=1}^n \mathbf{x}_j(k+1) \right), \label{eq:complete_dual_ascent_distributed_ascent}
\end{align}
\end{subequations}
which is an equivalent reformulation of \eqref{eq:equivalent_full}. 
Clearly, \eqref{eq:complete_dual_ascent_distributed_ascent} needs to be modified in order to account for a time varying (and not fully connected) network. For this purpose, we require a scheme that can track $\frac{1}{n}\sum_{j=1}^n\mathbf{x}_j(k+1)$ and that is distributed in nature. Hence, we have the following
\begin{subequations}
\label{eq:incomplete_dual_ascent}
\begin{align}
& \mathbf{x}(k+1) := & \text{arg }\underset{x \in \mathbb{R}^{np}}{\text{min}} & f(\mathbf{x}) - \mathbf{y}(k)^T\mathbf{x} \label{eq:incomplete_dual_ascent_obtain} \\
& \mathbf{z}(k+1) := && ? \\
& \mathbf{y}(k+1) := && \!\!\!\!\!\!\!\!\!\!\!\!\!\! \mathbf{y}(k) - c\left(\mathbf{x}(k+1) - \mathbf{z}(k+1)\right), \label{eq:incomplete_dual_ascent_ascent}
\end{align}
\end{subequations}
where we now need to establish which properties are desired of the iterate $\mathbf{z}(k).$ Clearly, we require that $\mathbf{z}(\infty) = \left(\boldsymbol{\Pi}_{\mathbf{1}_n} \otimes \mathbf{I}_p\right)\mathbf{x}(\infty).$ Further, we wish that $\mathbf{x}(\infty) = \mathbf{x}^{\star} = \left(\boldsymbol{\Pi}_{\mathbf{1}_n} \otimes \mathbf{I}_p \right)\mathbf{x}^{\star}.$ However, other properties may also be desirable. 
Let $\mathcal{Y} \triangleq \{\mathbf{y}: (\boldsymbol{\Pi}_{\mathbf{1}_n} \otimes \mathbf{I}_p)\mathbf{y} = \mathbf{0}\}.$ When dealing with the iterates in \eqref{eq:complete_dual_ascent_distributed} a bound on $\|\mathbf{x}(k+1) - \mathbf{x}(k)\|$ will imply a bound on the quantity $\|\mathbf{y}(k-1) -[\mathbf{y}(k-1) - c\mathbf{x}(k)]_{\mathcal{Y}} \|,$ where $[\cdot]_{\mathcal{Y}}$ denotes the projection on the set $\mathcal{Y},$ under suitable assumptions on $f.$ This is relevant because the quantity $\|\mathbf{y}(k-1) - [\mathbf{y}(k-1) - c \mathbf{x}(k)]_{\mathcal{Y}}\|$ is commonly used as a measure of sub-optimality of dual ascent \cite{luo} and can be used to bound $\|\mathbf{x}(k) - \mathbf{x}^{\star}\|.$ Hence, if $\mathbf{z}(k)$ is close to $\mathbf{z}(\infty)$ when the variations $\|\mathbf{x}(k+1) - \mathbf{x}(k)\|$ are small we are essentially preserving the relevance of bounding $\|\mathbf{x}(k+1) - \mathbf{x}(k)\|$ from the iterates in \eqref{eq:complete_dual_ascent_distributed} for the iterates in \eqref{eq:incomplete_dual_ascent}.

The work in \cite{interpret} uses the scheme of \cite{martinez} in order to track the average gradient of $f$ using the following scheme
\begin{align}
\hat{\nabla}f(\mathbf{x}(k+1)) = (\mathbf{W}(k)\otimes \mathbf{I}_p)\hat{\nabla}f(\mathbf{x}(k))  \nonumber\\ + \nabla f(\mathbf{x}(k+1)) - \nabla f(\mathbf{x}(k)),
\end{align}   
or alternatively
\begin{align}
\hat{\nabla}f_i(\mathbf{x}(k+1)) = \sum_{j \in \mathcal{N}_i(k) \cup \{i\}} w_{ij}(k)\hat{\nabla}f_j(\mathbf{x}_j(k)) \\
+ \nabla f_i(\mathbf{x}_i(k+1)) - \nabla f_i(\mathbf{x}_i(k)),
\end{align}
where $\mathbf{W}(k) \in \mathbb{R}^{n \times n}$ is a mixing matrix whose properties we will formally introduce later on. Element $(i,j)$ of $\mathbf{W}(k),$ i.e., $w_{ij}(k),$ is the weight node $i$ assigns to the information coming from node $j.$
This scheme has the two properties we desire, i.e. $\hat{\nabla}f(\mathbf{x}(\infty)) = (\boldsymbol{\Pi}_{\mathbf{1}_n} \otimes \mathbf{I}_p)\nabla f(\mathbf{x}(\infty))$ and as the differences $\nabla f(\mathbf{x}(k+1)) - \nabla f(\mathbf{x}(k)) $ become small $\hat{\nabla}f(\mathbf{x}(k+1))$ converges to $\hat{\nabla}f(\mathbf{x}(\infty))$ \cite{interpret}. We will therefore, for PANDA, use an analogous scheme to track $(\boldsymbol{\Pi}_{\mathbf{1}_n}\otimes \mathbf{I}_p)\mathbf{x}(k).$ The estimated network average at iterate $k+1$ will be
\begin{equation}
\mathbf{z}(k+1) := (\mathbf{W}(k) \otimes \mathbf{I}_p)\mathbf{z}(k) + \mathbf{x}(k+1) - \mathbf{x}(k),
\end{equation} 
which requires that the nodes exchange their previous estimate with that of their neighbor. The precise iterates with clarifications on what can be computed in a distributed manner and what information is to be exchanged are given in Algorithm \ref{alg:ddiging}, which
\floatstyle{spaceruled}
\restylefloat{algorithm}
\begin{algorithm}
\caption{PANDA \label{alg:ddiging}}
\begin{algorithmic}[1]
\State Choose step size $c > 0$ and pick $\mathbf{z}(0)= \mathbf{x}(0) = \mathbf{0}$ and $\mathbf{y}(0)$ such that $(\boldsymbol{\Pi}_{\mathbf{1}_n} \otimes \mathbf{I}_p)\mathbf{y}(0) = \mathbf{0}.$
\For{$k=0,1,\hdots$} Each agent $i$:
\State Computes $$
\mathbf{x}_i(k+1) := \text{arg } \underset{\mathbf{x}_i \in \mathbb{R}^{p}}{\text{min}}  f_i(\mathbf{x}_i) - \mathbf{y}_i(k)^T\mathbf{x}_i$$
\State Exchanges $\mathbf{z}_i(k)$
 with $\mathcal{N}_i(k).$ 
\State Computes 
$$\mathbf{z}_i(k+1) := \sum_{j \in \mathcal{N}_i \cup \{i\}}w_{ij}(k)\mathbf{z}_j(k) + \mathbf{x}_i(k+1) - \mathbf{x}_i(k)$$
\State Computes $$\mathbf{y}_i(k+1) := \mathbf{y}_i(k) - c(\mathbf{x}_i(k+1) -\mathbf{z}_i(k+1))$$
\EndFor
\end{algorithmic}
\end{algorithm}
 can be more compactly expressed as
 \begin{subequations}
 \label{eq:compact}
\begin{align}
&\mathbf{x}(k+1) := \text{arg } \underset{\mathbf{x} \in \mathbb{R}^{np}}{\text{min}} \quad f(\mathbf{x}) - \left(\mathbf{y}(k)\right)^T\mathbf{x} \label{eq:compact_primal}\\
& \!\mathbf{z}(k+1) := (\mathbf{W}(k) \otimes \mathbf{I}_p) \mathbf{z}(k) + \mathbf{x}(k+1) - \mathbf{x}(k) \label{eq:compact_averaging}\\
& \mathbf{y}(k+1) := \mathbf{y}(k) - c\big(\mathbf{x}(k+1) - \left(\mathbf{z}(k+1)\right)\big). \label{eq:iterates_bound}
\end{align}
\end{subequations}
Now that PANDA has been formally introduced we will provide its convergence properties. However, before doing so we will formally introduce the assumptions under which the convergence statement holds. This is done in the following section.
\section{Convergence Rate \label{section:convergence}}

In this section we formalize PANDA's convergence properties in the form of a theorem. For this we first provide the assumptions under which the theorem holds.
The assumptions present in this section are identical to the assumptions in \cite{interpret}.
To this end, consider a time-varying graph sequence $\mathcal{G}(k) = \{\mathcal{V},\mathcal{E}(k)\}.$ Note that the set of agents $\mathcal{V} = \{1,\hdots,n\}$ remains static. If node $i$ can communicate with node $j$ at time $k,$ $(i,j) \in \mathcal{E}(k).$ Let $\{\mathbf{W}(k)\}$ denote a mixing matrix sequence. 
\begin{assumption}[Mixing matrix sequence $\{\mathbf{W}(k)\}$ \cite{interpret}] \label{assumption:mixing} For any $k=0,1,\hdots,$ the mixing matrix $\mathbf{W}(k) \in \mathbb{R}^{n \times n}$ satisfies the following relations:
\begin{enumerate}
\item [(P1)] Decentralized property: if $i \neq j$ and $(j,i) \not \in \mathcal{E}(k)$ $w_{ji}(k) = 0.$
\item [(P2)] Double stochasticity: $\mathbf{W}(k)\mathbf{1}_n = \mathbf{1}_n,$ $\mathbf{1}_n^T\mathbf{W}(k) = \mathbf{1}_n^T.$
\item [(P3)] Joint spectrum property: Let  $\sigma_{\text{max}}(\cdot)$ denote the largest singular value of a matrix and let
\begin{equation}
\mathbf{W}_b(k) \triangleq \mathbf{W}(k)\mathbf{W}(k-1)\hdots\mathbf{W}(k-b+1)\, k\geq 0,
\end{equation}
for $k \geq 0$ and $b \geq k -1,$ with $\mathbf{W}_b(k) = \mathbf{I}_n$ for $k < 0$ and $\mathbf{W}_0(k) = \mathbf{I}_n.$ Then, there exists a positive integer $B$ such that
\begin{equation}
\label{eq:assumption:spectrum}
\underset{k \geq B-1}{\text{sup}} \sigma_{\text{max}} \left\{\mathbf{W}_B(k) - \frac{1}{n}\mathbf{1}_n\mathbf{1}_n^T\right\} < 1.
\end{equation}
\end{enumerate}
\end{assumption}
Realistic scenarios in which Assumption \ref{assumption:mixing} holds are provided in \cite{interpret}. Note however, that the assumption is much weaker than each $\mathcal{G}(k)$ being connected.

Note that as opposed to the communication matrices in Assumption \ref{assumption:communication} we do not require that the matrices $\mathbf{W}(k)$ are positive semi-definite. However, we do obtain better bounds for the convergence rate if the mixing matrices are picked to be positive semi-definite.

We will now introduce an assumption on $f$ which will in turn give us desired properties of $f^{*}.$ These will be announced and exploited in Section \ref{section:sketch}.
\begin{assumption}[Strong convexity and smoothness] \label{assumption:standard_linear}
The functions $f$ (cf. \eqref{eq:standard_dual_reformulation}) is proper closed strongly convex and Lipschitz differentiable, i.e. 
\begin{equation}
f(\mathbf{x}) \geq f(\mathbf{y}) + \left( \nabla f(\mathbf{y}) \right)^T(\mathbf{x} - \mathbf{y}) + \frac{\mu}{2}\|\mathbf{x} - \mathbf{y}\|^2
\end{equation}
and 
\begin{equation}
\|\nabla f(\mathbf{x}) - \nabla f(\mathbf{y})\|_2 \leq L \|\mathbf{x} - \mathbf{y}\|_2
\end{equation}
$\forall \mathbf{x},\mathbf{y} \in \mathbb{R}^p$ where $\mu > 0$ and $L < \infty$ are the strong convexity and Lipschitz constant respectively.
\end{assumption}
Assumption \ref{assumption:standard_linear} is a standard assumption when obtaining linear convergence using first order methods.

\begin{theorem}[PANDA converges R-linearly \label{theorem:linear}]
Let Assumptions \ref{assumption:mixing} and \ref{assumption:standard_linear} hold. Also, let
\begin{equation}
\delta \triangleq \underset{k \geq B-1}{\text{sup}} \,\left\{\sigma_{\text{max}} \left\{\mathbf{W}_B(k) - \frac{1}{n}\mathbf{1}_{n}\mathbf{1}_n^T \right\} \right\}.
\end{equation}
Finally let $\kappa \triangleq \frac{\mu}{L}$ denote the condition number of $f.$ Then, for any step-size
\begin{equation}
c \in \left(0,\frac{\mu \sqrt{\kappa}}{4B^2}(1-\delta^2) \right],
\end{equation}
the sequence $\{\mathbf{y}(k)\}$ converges to $\mathbf{y}^{\star},$ the unique solution of \eqref{eq:dual_fully_connected} and $\{\mathbf{x}(k)\}$ converges to $\mathbf{x}^{\star},$ the unique solution of \eqref{eq:complete_dual_reformulation} at a global R-linear rate $\mathcal{O}(\lambda^k),$ where $\lambda$ is given by
\begin{equation}
\label{eq:rates}
\lambda = 
\begin{cases}
\sqrt[2B]{1 - \frac{c}{2L}} & \text{if } c \in (0,\alpha] \\
\sqrt[B]{\delta + \sqrt{\frac{4cB^2}{\mu \sqrt{\kappa}}}} & \text{if } c \in \left( \alpha, \frac{\mu \sqrt{\kappa}}{4B^2}(1-\delta^2)\right],
\end{cases}
\end{equation} 
where 
\begin{equation}
\alpha \triangleq 2 \sqrt{\kappa} \mu \left( \frac{\sqrt{(1-\delta^2)\kappa^{2/3} + 8 B^2} - 8 \delta B}{\kappa^{3/2} + 8B^2} \right)^2. 
\end{equation}
\end{theorem}
A sketch of the proof of Theorem \ref{theorem:linear} is provided in Section \ref{section:sketch}. The full details are omitted due to space constraints.

Note that neither the convergence rates nor the step-size of PANDA depend explicitly on the size of the network. However, via the parameters $\delta$ and $B$ the rate and the step-size may depend on the network size. The bounds in \eqref{eq:rates} are not tight but they provide  a rate dependence on the objective function's condition number and network parameters. While it would be interesting to establish for which problem and network parameters PANDA obtains better convergence bounds than DIGing the amount and discrepancy of parameters in \eqref{eq:rates} and \cite{interpret}, and the additional tunable parameters in \cite{interpret} make the task of obtaining an explicit expression difficult. This is therefore left for future work. 

\section{Proof Sketch \label{section:sketch}}
In this section we provide a sketch of the proof of Theorem \ref{theorem:linear}. This proof essentially relies on the small gain theorem and has many elements in common with the convergence proof in \cite{interpret}. Due to space restrictions we do not provide the entire proof here but do discuss the main steps involved and provide some intuition behind it. 

Just like \cite{interpret} we require the small gain theorem, which we introduce here for completeness. However, we first need to define two norms. Let $\mathbf{s}^i  = \{\mathbf{s}^{i}(0),\mathbf{s}^{i}(1),\hdots\}$ denote an infinite sequence of vectors $\mathbf{s}^i(k) \in \mathbb{R}^{np},$ $\forall i.$ Further, let
\begin{equation}
\|\mathbf{s}^i\|^{\lambda,K} \triangleq \underset{k=0,\hdots,K}{\text{sup}} \quad \frac{1}{\lambda^k} \|\mathbf{s}^{i}(k)\|
\end{equation}
and 
\begin{equation}
\|\mathbf{s}^i\|^{\lambda} \triangleq \underset{k \geq 0}{\text{sup}}\frac{1}{\lambda^k}\|\mathbf{s}^{i}(k)\|.
\end{equation}
\begin{theorem}[Small Gain Theorem \cite{interpret}]
Suppose $\mathbf{s}^1,\hdots,\mathbf{s}^m$ are vector sequences such that for all positive integers $K$ and for each $i = 1,\hdots,m,$ we have an arrow $\mathbf{s}^i \rightarrow \mathbf{s}^{(i \!\!\! \mod m) + 1},$ that is,
\begin{equation}
\label{eq:cyclic}
\|\mathbf{s}^{(i \!\!\!\!\!\mod m) + 1}\|^{\lambda,K} \leq \gamma_i \|\mathbf{s}^i\|^{\lambda,K} + \omega_i,
\end{equation}
where the constants $\gamma_1,\hdots,\gamma_m$ and $\omega_1,\hdots,\omega_m$ are independent of $K$. Further, suppose that the constants $\gamma_1,\hdots,\gamma_m$ are nonnegative and satisfy
\begin{equation}
\label{eq:small_gain}
\gamma_1\gamma_2\hdots\gamma_m < 1.
\end{equation}
Then we have that
\begin{align}
&\|\mathbf{s}^1\|^{\lambda} \leq \frac{1}{1-\gamma_1\gamma_2 \hdots \gamma_m}\big(\omega_1\gamma_2\gamma_3\hdots\gamma_m + \omega_2\gamma_3\gamma_4\hdots\gamma_m  \nonumber\\
 &+ \hdots + \omega_{m-1}\gamma_m + \omega_m \big).
\end{align}
Due to the cyclic nature of \eqref{eq:cyclic} similar bounds hold for the remaining sequences.
\end{theorem}
Note that if a sequence $\mathbf{s}^i$ fulfills $\|\mathbf{s}^i\|^{\lambda}  < \infty,$ the sequence converges to $\mathbf{0}$ geometrically fast at rate $\mathcal{O}(\lambda^k)$ \cite{interpret}.
By establishing \eqref{eq:cyclic} and \eqref{eq:small_gain} for every $K$ and some $\lambda < 1$ the small gain theorem can be used to conclude that the sequences $\mathbf{s}^1,\hdots,\mathbf{s}^m$ converge to $\mathbf{0}$ at rate $\mathcal{O}(\lambda^k).$

In the context of PANDA and DIGing, one must define sequences of residuals such that we can derive a bound of the type of \eqref{eq:small_gain} for every $K.$ Then, if the gains fulfill \eqref{eq:small_gain} we will be capable of establishing the geometric decrease of all residuals. Obviously, an essential residual of interest is the dual residual, i.e. $\mathbf{y}(k) - \mathbf{y}^{\star}.$  The sequences of residuals we have used to establish Theorem \ref{theorem:linear} are
\begin{subequations}
\begin{align}
\mathbf{r}(k) \triangleq \mathbf{y}(k) - \mathbf{y}^{\star} \\
\mathbf{x}^{\perp}(k+1) \triangleq (\boldsymbol{\Pi}_{\mathbf{1}_n}^{\perp}\otimes \mathbf{I}_p) \mathbf{x}(k+1) \\
\Delta \mathbf{y}(k+1) \triangleq \mathbf{y}(k+1) - \mathbf{y}(k) \\
\mathbf{z}^{\perp}(k+1) \triangleq (\boldsymbol{\Pi}_{\mathbf{1}_n}^{\perp} \otimes \mathbf{I}_p) \mathbf{z}(k+1) \\
\boldsymbol{\Delta}_{xz}^{\perp}(k+1) \triangleq \mathbf{x}^{\perp}(k+1) - \mathbf{z}^{\perp}(k+1),
\end{align}
\end{subequations}
for $k \geq 0.$ We adopt the convention $\mathbf{x}^{\perp}(0) = \Delta \mathbf{y}(0) = \mathbf{z}^{\perp}(0) = \boldsymbol{\Delta}_{xz}^{\perp} = \mathbf{0}.$
Since the goal is to ultimately establish that $\|\mathbf{r}(k)\| \to 0$ R-linearly we will start the circle of arrows with the sequence $\mathbf{r}(k)$ and then proceed as
\begin{equation}
\mathbf{r} \rightarrow \mathbf{x}^{\perp} \rightarrow \boldsymbol{\Delta}^{\perp}_{xz}  \rightarrow \Delta\mathbf{y} \rightarrow  \mathbf{z}^{\perp} \rightarrow \mathbf{r}.
\end{equation}
More specifically we establish the following relations
\begin{enumerate}
\item [(A1)] $\|\mathbf{x}^{\perp}\|^{\lambda,K} \leq \gamma_1 \|\mathbf{r}\|^{\lambda,K} + \omega_1,$ where $\gamma_1 = \frac{1}{\mu \lambda}$ and $\omega_1 =0.$
\item [(A2)] $\|\boldsymbol{\Delta}_{xz}^{\perp}\|^{\lambda,K} \leq \gamma_2 \|\mathbf{x}^{\perp}\|^{\lambda,K} + \omega_2,$ where $\gamma_2 = \frac{2(1-\lambda^B)}{(1-\lambda)(\lambda^B - \delta)}$ and $\omega_2 = \frac{\lambda^B}{\lambda^B - \delta}\sum_{t=1}^B \lambda^{1-t}\|\mathbf{x}^{\perp}(t-1)\|,$ for $\lambda^B > \delta.$
\item [(A3)] $\|\Delta \mathbf{y}\|^{\lambda,K} \leq \gamma_3 \|\boldsymbol{\Delta}_{xz}^{\perp}\|^{\lambda,K} + \omega_3, $ where $\gamma_3 = c$ and $\omega_3 = 0.$
\item [(A4)] $\|\mathbf{z}^{\perp}\|^{\lambda,K} \leq \gamma_4 \|\Delta \mathbf{y}\|^{\lambda,K} + \omega_4,$ where $\gamma_4 = \frac{(1-\lambda^B)}{\mu(1-\lambda)(\lambda^B - \delta)}$ and $\omega_4 = \frac{\lambda^B}{\lambda^B - \delta} \sum_{t=1}^B \lambda^{1-t}\|\mathbf{z}(t-1)\|$ for $\lambda > \delta.$
\item [(A5)] $\|\mathbf{r}\|^{\lambda,K} \leq \gamma_5 \|\mathbf{z}^{\perp}\| + \omega_5,$ where $\gamma = \sqrt{L\mu}$ and $\omega_5 = 2\|\mathbf{r}(0)\|$ for $\lambda \in \left[\sqrt{1-\frac{c}{2L}},1 \right)$ and $c \in \left(0,\frac{\mu}{2}\right].$ 
\end{enumerate}

Under Assumption \ref{assumption:standard_linear} the convex conjugate $f^*$ can be shown to be $\frac{1}{L}-$strongly convex and $\frac{1}{\mu}-$Lipschitz continuous. Relation (A1) can be shown to hold true by relying on this fact. This establishes that we can use the dual residual to bound the disagreement within the network. 

 Relation (A2) relies on Assumption \ref{assumption:mixing} which gives the sequence $\{\mathbf{W}(k)\}$ some averaging properties. Using these properties and setting $\lambda^B > \delta$ the bound is obtained following a similar procedure to that in \cite{interpret}. (A3) simply uses the iterates \eqref{eq:iterates_bound} to bound the variation $\|\mathbf{y}(k) - \mathbf{y}(k-1)\|.$ Note that the inclusion of $c$ in the term $\gamma_1\hdots\gamma_5$ is critical to the fulfillment of \eqref{eq:small_gain}. (A2) and (A3) together imply that the more the nodes of the network agree the less two consecutive dual iterates differ from each other. This in turn will imply that the quantity $\|\mathbf{y}(k-1)-[\mathbf{y}(k-1)-c\mathbf{x}(k) +c\mathbf{z}(k)]_{\mathcal{Y}}\|$ is small. Recall from the discussion in Section \ref{section:algorithm} that in classic dual ascent the quantity $\|\mathbf{y}(k-1)-[\mathbf{y}(k-1)-c\mathbf{x}(k)]_{\mathcal{Y}}\|$ is used as a sub-optimality measure.  The two following arrows, (A4) and (A5), are therefore  aimed at relating PANDA to dual ascent. In order to do this we equivalently re-write the iterate in \eqref{eq:iterates_bound} as
\begin{equation}
\mathbf{y}(k+1) := \mathbf{y}(k) - c(\boldsymbol{\Pi}_{\mathbf{1}_n}^{\perp} \otimes \mathbf{I}_p)(\mathbf{x}(k+1) - \mathbf{z}^{\perp}(k+1)). \label{eq:gradient_error}
\end{equation} 
The proof of this statement is omitted due to space restrictions but it relies on the shape of the iterates $\{\mathbf{z}(k)\},$ the properties of the sequence $\{\mathbf{W}(k)\}$ and the fact that $\boldsymbol{\Pi}_{\mathbf{1}_n}^{\perp}$ is idempotent. 

(A4) establishes that the \emph{gradient error} $\mathbf{z}^{\perp}(k+1)$ is small if the quantity $\|\Delta \mathbf{y}(k+1)\|$ is small. Note that by seeing \eqref{eq:gradient_error} as dual ascent with an oracle that provides inexact gradient information a bound on $\|\mathbf{z}^{\perp}(k+1)\|$ implies a bound on the error. This interpretation is exploited by using the framework in \cite{oracle} to establish the bound (A5).  Finally, we find values of $\lambda$ and $c$ that fulfill all the restrictions in (A1)-(A5) and that lead to $\gamma_1\hdots\gamma_5 < 1$ concluding this proof sketch.
\section{Numerical Experiments \label{section:numerical}}
In this section we experimentally verify the theoretical findings in the previous sections. We consider a network in which agent $i \in \{1,\hdots,n\}$ has its own measurement $\mathbf{m}^{(i)},$ measurement equation
\begin{equation}
\mathbf{m}^{(i)} = \mathbf{H}^{(i)}\mathbf{x} + \mathbf{e}^{(i)},
\end{equation} 
where $\mathbf{H}^{(i)} \in \mathbb{R}^{p\times p}$ denotes agent $i'$s measurement matrix and $\mathbf{e}^{(i)}$ denotes its measurement noise. The goal is for the nodes to cooperatively estimate $\mathbf{x}.$ This will be done by solving the optimization problem
\begin{align}
\label{eq:numerical_optimization}
\underset{\mathbf{x} \in \mathbb{R}^{p}}{\text{minx}} \quad \sum_{i=1}^n \frac{1}{2}\|\mathbf{H}^{(i)}\mathbf{x} -\mathbf{m}^{(i)}\|^2.
\end{align}
The matrices $\mathbf{H}^{(i)}$ are set to be of size $5 \times 5$ and have been, just as $\mathbf{w}^{(i)},$  randomly generated following the standard Gaussian distribution. We have generated as many random matrices required to obtain all $\mathbf{H}^{(i)}$ with a condition number smaller than or equal to 100. 

The sequence of graphs $\mathcal{G}(k)$ is modeled using an i.i.d. stochastic process. Each $\mathcal{G}(k)$ is generated based on the fully connected network where each of the links will be independently removed with probability 0.2.  
\begin{figure}
\vspace{10pt}
\scalebox{0.57}{\input{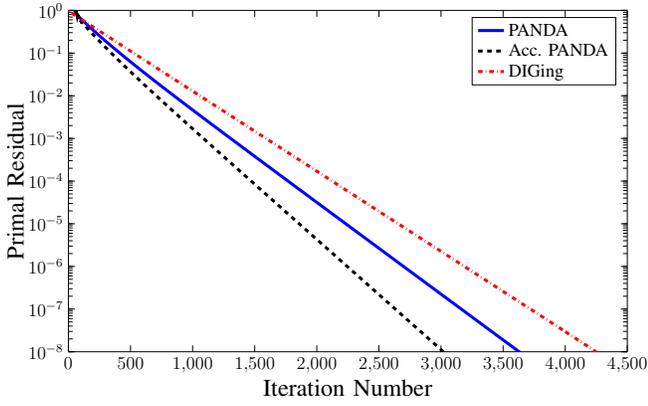}}
\caption{Plot of residuals $\frac{\|\mathbf{x}(k) - \mathbf{x}^{\star}\|}{\|\mathbf{x}^{\star} - \mathbf{x}(0)\|}$ for a time varying graph. The step-sizes are $c = 0.013$ for PANDA and its accelerated version and $\alpha = 0.24$ for DIGing. The acceleration weight for accelerated PANDA is chosen to be 0.2. \label{fig:plot}}
\end{figure}
For a randomly generated optimization problem \eqref{eq:numerical_optimization} the step-sizes and acceleration weight for PANDA and DIGing have been \emph{hand optimized} in order to obtain the plots in Figurcontrol and decision conferencee \ref{fig:plot}. From Figure \ref{fig:plot} we can see that PANDA converges at a faster rate than DIGing does. We can also see that we do obtain a gain in the speed of convergence by using Nesterov's acceleration scheme \cite{nesterov,optimal}. In particular, we have used the scheme
\begin{subequations}
\begin{align}
\mathbf{x}(k+1) := \text{arg }\underset{\mathbf{x} \in \mathbb{R}^{np}}{\text{min}} \quad f(\mathbf{x}) - \mathbf{\bar{y}}(k)^T\mathbf{x} \\
\mathbf{z}(k+1) := (\mathbf{W}(k) \otimes \mathbf{I}_p) \mathbf{z}(k) + \mathbf{x}(k+1) - \mathbf{x}(k) \\
\mathbf{y}(k+1) := \mathbf{y}(k) - c (\mathbf{x}(k+1) - \mathbf{z}(k+1))\\
\mathbf{\bar{y}}(k+1) := (1 + \eta)\mathbf{y}(k+1) - \eta \mathbf{y}(k),
\end{align}
\end{subequations}
for $\eta = 0.2.$
When it comes to communication cost, recall that both PANDA and its accelerated counterpart require the exchange of half as many real numbers per iteration as DIGing implying that the communication cost per iterate is cut by half. Since PANDA does not converge slower than DIGing this implies that the only downside of using PANDA is the additional computational cost each of the iterates have. We can conclude therefore that we are trading communication for computation.
\section{Conclusions and Further Work \label{section:conclusions}} 
In this paper we proposed PANDA, a dual ascent based method for time varying graphs. We studied its convergence properties and compared its experimental performance to that of DIGing.
 One of the advantages of PANDA is that it requires communicating half as many quantities as DIGing per iteration. On the other hand, PANDA's iterates are in general computationally more expensive than those of DIGing. This implies that while DIGing may be more suitable for scenarios in which the computation cost is high PANDA may be more suitable for scenarios in which the communication cost is high. We have also experimentally seen that we can speed up PANDA's convergence by using Nesterov's acceleration scheme. Proving this is left for future work.
\input{apendix}

\end{document}

%% file: apendix.tex
\begin{appendix}
This appendix is attached for review purposes only and, if accepted, will be removed from the paper. We are ware that the sketches of the proofs contained here contain only the essential steps to establish the results and are not detailed enough to be followed comfortably. We will use that $f^*$ is $\frac{1}{L}-$strongly convex and has $\frac{1}{\mu}-$Lipschitz continuous gradients, which can be proven by combining Theorem 6 in \cite{dual_smooth} and Theorem 2.1.5. in \cite{nesterov}. We also require the use of Lemma 3.4 in \cite{interpret} to establish both (A2) and (A4).
\section*{Proof sketch of (A1)}
Since $f^*$ has $\frac{1}{\mu}-$Lipschitz continuous gradients we have
\begin{equation}
\label{eq:LHS}
\|\mathbf{x}(k+1) - \mathbf{x}^{\star}\| \leq \frac{1}{\mu}\|\mathbf{y}(k) - \mathbf{y}^{\star}\|,
\end{equation}
since $((\boldsymbol{\Pi}_{\mathbf{1}_n}^{\perp}\otimes \mathbf{I}_p)\mathbf{x}(k+1))^T (\mathbf{x}^{\star} - (\boldsymbol{\Pi}_{\mathbf{1}_n}\otimes \mathbf{I}_p)\mathbf{x}(k+1)) = \mathbf{0},$ by adding and subtracting $\boldsymbol{\Pi}_{\mathbf{1}_n}\mathbf{x}(k+1)$ in the LHS of \eqref{eq:LHS} it can be established that 
\begin{equation}
\|\mathbf{x}^{\perp}(k+1)\| \leq \frac{1}{\mu}\|\mathbf{r}(k)\|.
\end{equation}
By following a procedure similar to that of the proof of Lemma 3.9 in \cite{interpret} we establish
\begin{equation}
\|\mathbf{x}^{\perp}\|^{\lambda,K} \leq \frac{1}{\lambda \mu}\|\mathbf{r}\|^{\lambda,K}, \, \forall \lambda \in (0,1).
\end{equation}
\section*{Proof sketch of (A2)}
By recursively applying the iterate \eqref{eq:compact_averaging} on the definition of $\boldsymbol{\Delta}_{xz}^{\perp}(k+1),$ we have that
\begin{align}
&\boldsymbol{\Delta}_{xz}^{\perp}(k+1) = (\boldsymbol{\Pi}_{\mathbf{1}_n}^{\perp}\mathbf{W}_B(k) \otimes \mathbf{I}_p) \boldsymbol{\Delta}_{xz}^{\perp}(k+1-B) + \nonumber\\
& \sum_{t=k+1-B}^k (\boldsymbol{\Pi}_{\mathbf{1}_n}^{\perp}\mathbf{W}_{k-(t+1)}(k)(\mathbf{I}_n - \mathbf{W}(t))\otimes \mathbf{I}_p)\mathbf{x}^{\perp}(t).
\end{align}
By using the triangle inequality and the fact that $\{\mathbf{W}(k)\}$ is doubly stochastic we obtain the following bound
\begin{align}
&\|\boldsymbol{\Delta}_{xz}^{\perp}(k+1)\| \leq \|(\mathbf{W}_B(k) - \boldsymbol{\Pi}_{\mathbf{1}_n}\otimes \mathbf{I}_p) \boldsymbol{\Delta}_{xz}^{\perp}(k+1-B)\| \nonumber\\
&+ \sum_{t=k+1-B}^k \|\boldsymbol{\Pi}_{\mathbf{1}_n}^{\perp}(\mathbf{W}_{k-(t+1)}(k)(\mathbf{I}_n - \mathbf{W}(t))\otimes \mathbf{I}_p)\mathbf{x}^{\perp}(t)\|. \label{eq:bounded_term}
\end{align}
Again since $\{\mathbf{W}(k)\}$ are doubly stochastic it holds that $\|(\boldsymbol{\Pi}_{\mathbf{1}_n}^{\perp}\mathbf{W}_{k-(t+1)}(k)(\mathbf{I}_n - \mathbf{W}(t)) \otimes \mathbf{I}_p)\mathbf{x}^{\perp}(t)\| \leq \rho (\mathbf{W}_{k-(t+1)}(k))\rho(\mathbf{I}_n - \mathbf{W}(t))\|\mathbf{x}^{\perp}(t)\| \leq 2\|\mathbf{x}^{\perp}(t)\|,$ where $\rho(\cdot)$ denotes the spectral radius. Hence, \eqref{eq:bounded_term} is upper bounded by
\begin{equation}
\label{eq:bound_on_term}
\sum_{t=k+1-B}^k 2\|\mathbf{x}^{\perp}(t)\|.
\end{equation}
Note that if the sequence of matrices $\{\mathbf{W}(k)\}$ is such that each matrix is positive semi-definite the factor 2 in \eqref{eq:bound_on_term} is replaced by a 1 since $\rho(\mathbf{I}_n - \mathbf{W}(t)) \leq 1$.
From here a procedure similar to that in the proof of Lemma 3.10. in \cite{interpret} can be followed to yield
\begin{align}
\|\boldsymbol{\Delta}_{xz}^{\perp}\|^{\lambda,K} \leq \frac{2(1-\lambda^B)}{(1-\lambda)(\lambda^B - \delta)}\|\mathbf{x}^{\perp}\|^{\lambda,K} + \\
\frac{\lambda^B}{\lambda^B - \delta}\sum_{t=1}^B \lambda^{1-t}\|\mathbf{x}^{\perp}(t-1)\|, \, \forall \lambda \in (\delta^{1/B},1).\nonumber
\end{align}
Note that the term $\sum_{t=1}^B\lambda^{1-t}\|\mathbf{x}^{\perp}(t-1)\|$ is bounded since it involves a bounded number of bounded iterates.
\section*{Proof sketch of (A3)}
(A3) follows from applying iterate \eqref{eq:iterates_bound}.
\section*{Proof sketch of (A4)}
Using the iterates \eqref{eq:compact_averaging} on the definition of $\mathbf{z}^{\perp}(k+1),$ $\mathbf{z}^{\perp}(k+1)$ can be expressed as 
\begin{align}
\mathbf{z}^{\perp}(k+1) = (\boldsymbol{\Pi}_{\mathbf{1}_n}^{\perp}\mathbf{W}_b(k)\otimes \mathbf{I}_p) \mathbf{z}^{\perp}(k+1-B)  \\
+ \sum_{t=k-B+1}^{k} (\boldsymbol{\Pi}_{\mathbf{1}_n}^{\perp}\mathbf{W}_{k-t}(k) \otimes \mathbf{I}_p)(\mathbf{x}(t+1) - \mathbf{x}(t)). \nonumber
\end{align}
Following similar steps as in the Proof sketch of (A2) we obtain
\begin{align}
\|\mathbf{z}^{\perp}(k+1)\| \leq \delta\|\mathbf{z}^{\perp}(k+1-B)\| + \\
\sum_{t=1}^B \|\mathbf{x}(k+2-t) - \mathbf{x}(k+1-t)\|. \nonumber
\end{align}
Since $f^{*}$ has $\frac{1}{\mu}-$Lipschitz continuous gradients we have that
\begin{equation}
\|\mathbf{x}(k+2-t) - \mathbf{x}(k+1-t)\| \leq \frac{1}{\mu}\|\mathbf{y}(k+1-t) - \mathbf{y}(k-t)\|,
\end{equation}
implying, following the procedure from the Proof sketch of (A2), that
\begin{align}
\|\mathbf{z}^{\perp}\|^{\lambda,K} \leq \frac{(1-\lambda^B)}{\mu(1-\lambda)(\lambda^B - \delta)}\|\Delta \mathbf{y}\|^{\lambda,K} \\
+ \frac{\lambda^B}{\lambda^B - \delta}\sum_{t=1}^B \lambda^{1-t}\|\mathbf{z}^{\perp}(t-1)\|, \,\forall\lambda \in (\delta^{1/B},1)\nonumber
\end{align}
\section*{Proof sketch of (A5)}
Before proving that (A5) is true we establish that the iterate \eqref{eq:iterates_bound} can indeed be expressed as \eqref{eq:gradient_error}. To see this Let us express $\mathbf{x}(k+1) - \mathbf{z}(k+1)$ as a function of the iterates $\{\mathbf{x}(t)\}_{t=0}^{k}.$ In particular, by applying \eqref{eq:compact_averaging} iteratively we obtain
\begin{align}
&\mathbf{x}(k+1) - \mathbf{z}(k+1) = ((\mathbf{I}_n - \mathbf{W}(k)) \otimes \mathbf{I}_p)\mathbf{x}(k) + \nonumber \\
&(\mathbf{W}(k)(\mathbf{I}_n - \mathbf{W}(k-1)) \otimes \mathbf{I}_p)\mathbf{x}(k-1) +  \\
&\sum_{t=0}^{k-3}(\mathbf{W}(k)\hdots\mathbf{W}(t+2)(\mathbf{I}_n - \mathbf{W}(t+1)) \otimes \mathbf{I}_p)\mathbf{x}(t+1), \nonumber
\end{align}
where we have only used the fact that $(\mathbf{A}\otimes \mathbf{B})(\mathbf{C} \otimes \mathbf{D}) = \mathbf{AC}\otimes\mathbf{BD}$ for matrices of compatible dimensions. Then, the equivalence between \eqref{eq:iterates_bound} and \eqref{eq:gradient_error} follows from the sequence of matrices $\{\mathbf{W}(k)\}$ being doubly stochastic and $\boldsymbol{\Pi}_{\mathbf{1}_n}^{\perp}$ being idempotent.
Further, \eqref{eq:gradient_error} can be equivalently written as
\begin{align}
\label{eq:inexact_oracle}
\mathbf{y}(k+1) := \text{arg }\underset{\mathbf{y} \in \mathcal{Y}}{\text{min}} \quad \left(\nabla f^*(\mathbf{y}(k)) + \boldsymbol{\epsilon}(k) \right)^T(\mathbf{y} - \mathbf{y}(k)) + \nonumber \\
 \frac{1}{2c}\|\mathbf{y} - \mathbf{y}(k)\|^2,
\end{align}
where $\boldsymbol{\epsilon}(k)$ models the error in $\nabla f^*(\mathbf{y}(k)).$ Note that $\boldsymbol{\epsilon}(k) = \mathbf{z}^{\perp}(k+1).$ 
The iterate in \eqref{eq:inexact_oracle} is an iterate in the same form as that provided in \cite{oracle} to analyze the implications of an inexact oracle in projected gradient descent for smooth and strongly convex objective functions. However, the analysis is not identical since the error term $\boldsymbol{\epsilon}(k)$ does not have a universal upper bound as opposed to the error terms in \cite{oracle}. As in \cite{oracle} we will first establish bounds that are a direct consequence of the Lipschitz continuity of the gradients of $f^*$ and the strong convexity of $f^*$ that are tweaked to take into account the gradient error $\boldsymbol{\epsilon}(k).$ In particular, we have that since $f^*$ has $\frac{1}{\mu}-$Lipschitz continuous gradients
\begin{align}
&f^*(\mathbf{y}) \leq f^*(\mathbf{y}(k)) - \boldsymbol{\epsilon}(k)^T(\mathbf{y}-\mathbf{y}(k)) \\
&+ (\nabla f^{*}(\mathbf{y}(k)) + \boldsymbol{\epsilon}(k))^T(\mathbf{y} - \mathbf{y}(k)) + \frac{1}{2\mu}\|\mathbf{y}(k) - \mathbf{y}\|^2 \leq \nonumber \\
&f^*(\mathbf{y}(k)) + (\nabla f^*(\mathbf{y}(k)) + \boldsymbol{\epsilon}(k))^T(\mathbf{y}(k+1) - \mathbf{y}(k)) \nonumber \\
&+ \|\boldsymbol{\epsilon}(k)\|\|\mathbf{y}(k) - \mathbf{y}\| + \frac{1}{2\mu}\|\mathbf{y} - \mathbf{y}(k)\|^2, \forall \mathbf{y} \in \mathcal{Y}. \nonumber
\end{align}
Using the Peter-Paul inequality $\|\boldsymbol{\epsilon}(k)\|\|\mathbf{y}(k) - \mathbf{y}\| \leq \frac{\mu \|\boldsymbol{\epsilon}(k)\|}{2} + \frac{1}{2\mu}\|\mathbf{y} - \mathbf{y}(k)\|^2$ yields
\begin{align}
\label{eq:Lipschitz_bound}
& f^*(\mathbf{y}) - f^*(\mathbf{y}(k)) \leq \\
&(\nabla f^*(\mathbf{y}(k)) + \boldsymbol{\epsilon}(k))^T(\mathbf{y} - \mathbf{y}(k)) + \nonumber \\
&\frac{1}{\mu}\|\mathbf{y}(k) - \mathbf{y}\|^2 + \frac{\mu \|\boldsymbol{\epsilon}(k)\|^2}{2},\, \forall \mathbf{y} \in \mathcal{Y}. \nonumber
\end{align}
Analogously by using the fact that $f^*$ is $\frac{1}{L}-$strongly convex we have
\begin{align}
\label{eq:strong_bound}
f^*(\mathbf{y}) - f^*(\mathbf{y}(k)) \geq (\nabla f^*(\mathbf{y}(k)) + \boldsymbol{\epsilon}(k))^T(\mathbf{y} - \mathbf{y}(k)) + \nonumber\\
\frac{1}{4L}\|\mathbf{y} - \mathbf{y}(k)\|^2 - L\|\boldsymbol{\epsilon}(k)\|^2, \forall \mathbf{y} \in \mathcal{Y}.
\end{align}
Combining both \eqref{eq:Lipschitz_bound} and \eqref{eq:strong_bound} yields
\begin{align}
\label{eq:all_bound}
\frac{1}{4L}\|\mathbf{y} - \mathbf{y}(k)\| \leq f^*(\mathbf{y}) - f_{\text{err}}^*(\mathbf{y}(k),\boldsymbol{\epsilon}(k)) \\ - (\nabla f^*(\mathbf{y}(k)) + \boldsymbol{\epsilon}(k))^T(\mathbf{y} - \mathbf{y}(k)) \nonumber \\
\leq \frac{1}{\mu}\|\mathbf{y}(k) - \mathbf{y}\|^2 + \frac{\mu \|\boldsymbol{\epsilon}(k)\|^2}{2},\,\forall \mathbf{y} \in \mathcal{Y}, \nonumber
\end{align}
where $f_{\text{err}}^*(\mathbf{y},\boldsymbol{\epsilon}) \triangleq f^*(\mathbf{y}) - L \|\boldsymbol{\epsilon}(k)\|^2.$ Now that we have established \eqref{eq:all_bound} we will proceed to analyze the iterates \eqref{eq:inexact_oracle}. For notational convenience let $r^2(k) \triangleq \|\mathbf{r}(k)\|^2.$ Then, for $k \geq 0$ it holds that
\begin{align}
\label{eq:triangle}
& r^2(k+1) = \|\mathbf{y}(k+1) - \mathbf{y}^{\star}\|^2 = r^2(k) \\
& + 2(\mathbf{y}(k+1) - \mathbf{y}(k))^T(\mathbf{y}(k+1) - \mathbf{y}^{\star})\nonumber\\
& - \|\mathbf{y}(k+1) - \mathbf{y}(k)\|^2. \nonumber
\end{align}
By writing the optimality condition to \eqref{eq:inexact_oracle} we obtain
\begin{align}
&(\nabla f^*(\mathbf{y}(k)) + \boldsymbol{\epsilon}(k) + \frac{1}{c}(\mathbf{y}(k+1) - \mathbf{y}(k))^T \\
& (\mathbf{y} - \mathbf{y}(k+1)) \geq 0,  \, \forall \mathbf{y} \in \mathcal{Y}, \nonumber
\end{align}
which can be particularized for $\mathbf{y}^{\star}$ yielding
\begin{align}
\label{eq:opt_condition}
(\mathbf{y}(k+1) - \mathbf{y}(k))^T(\mathbf{y}(k+1) - \mathbf{y}^{\star}) \\
\leq c (\nabla f^*(\mathbf{y}(k)) + \boldsymbol{\epsilon}(k))^T(\mathbf{y}^{\star} - \mathbf{y}(k+1)).
\end{align}
Using \eqref{eq:triangle} and \eqref{eq:opt_condition} yields
\begin{align}
\label{eq:ready_combine}
&r^2(k+1) \leq r^2(k) + 2c (\nabla f^*(\mathbf{y}(k)) + \boldsymbol{\epsilon}(k))^T(\mathbf{y}^{\star} - \mathbf{y}(k)) \nonumber\\
&- 2c(\nabla f^*(\mathbf{y}(k)) - \boldsymbol{\epsilon}(k))^T(\mathbf{y}(k+1) - \mathbf{y}(k)) \nonumber \\
&+ \frac{1}{2c}\|\mathbf{y}(k+1) - \mathbf{y}(k)\|^2. 
\end{align}
By using \eqref{eq:all_bound} with $\mathbf{y} = \mathbf{y}(k+1)$ on \eqref{eq:ready_combine} it holds that
\begin{align}
&r^2(k+1) \leq r^2(k) + 2c(\nabla f^*(\mathbf{y}(k)) + \boldsymbol{\epsilon}(k))^T(\mathbf{y}^{\star} - \mathbf{y}(k)) \nonumber \\
&- 2c(f^*(\mathbf{y}(k+1)) - f_{\text{err}}^*(\mathbf{y}(k),\boldsymbol{\epsilon}(k)) - \frac{\mu}{2}\|\boldsymbol{\epsilon}(k)\|^2),
\end{align}
as long as $c \in \left(0,\frac{\mu}{2} \right].$ By using \eqref{eq:all_bound} with $\mathbf{y} = \mathbf{y}^{\star}$ we further obtain 
\begin{align}
&r^2(k+1) \leq (1 - \frac{c}{2L})r^2(k) + \\
&2c(f^*(\mathbf{y}^{\star}) - f^{*}(\mathbf{y}(k+1))) + c \mu \|\boldsymbol{\epsilon}(k)\|^2. \nonumber
\end{align}
From here using the fact that $\nabla f^*(\mathbf{y}^{\star})^T(\mathbf{y} - \mathbf{y}^{\star}) = 0,\,\forall \mathbf{y} \in \mathcal{Y}$ we proceed as in the proof of Lemma 3.12 in \cite{interpret} yielding
\begin{equation}
\|\mathbf{r}\|^{\lambda,K} \leq 2 \|\mathbf{r}(0)\| + \sqrt{L\mu} \|\mathbf{z}^{\perp}\|^{\lambda,K},
\end{equation}
for $\lambda \in \left[\sqrt{1-\frac{c}{2L}},1 \right)$ and $c \in (0,\frac{\mu}{2}].$
\section*{Proof sketch of Theorem \ref{theorem:linear}}
The proof of Theorem \ref{theorem:linear} relies on using the statements (A1)-(A5) and showing that one can select parameters $c$ and $\lambda$ such that \eqref{eq:small_gain} is fulfilled. In particular the following set of inequalities need to hold
\begin{subequations}
\label{eq:all_restrictions}
\begin{align}
\frac{c(1-\lambda^B)^2}{\lambda(1-\lambda)^2(\lambda^B - \delta)^2} \leq \frac{\mu}{2}\sqrt{\frac{\mu}{L}} \label{eq:decreasing}\\
0 < c \leq \frac{\mu}{2} \\
\lambda \geq \sqrt{1 - \frac{c}{2L}} \\
\delta^{1/B} < \lambda < 1.
\end{align}
\end{subequations}
Note that \eqref{eq:decreasing} is monotonically decreasing in $\lambda.$ Hence one can always select a sufficiently small $c$ and large $\lambda < 1$ to enforce that \eqref{eq:all_restrictions} hold. While this guarantees that PANDA converges R-linearly it does not provide information of how to select the step-size nor the convergence rate. In order to obtain an explicit expression we renounce some tightness on the bounds.
In particular, we use the fact that for $0.5 \leq \lambda \leq 1$ it holds that 
$$(1-\lambda^B)^2/(\lambda(1-\lambda)^2) \leq 2B^2.$$ From here the analysis is very similar to that in \cite{interpret}.
\end{appendix}

%% file: panda.bbl
\begin{thebibliography}{99}
\bibitem{xiao}
J.J~.Xiao, A.~Ribeiro, et. al. ``Distributed compression-estimation using wireless sensor networks," \emph{IEEE Signal Proc. Magazine,} vol. 23, no. 4, pp. 27-41, July 2006.
\bibitem{decentralized_gradient}
K.~Yuan, Q.~Ling, et. al. ``On the Convergence of Decentralized Gradient Descent," \emph{SIAM. J. Optim.,} vol. 26, no. 3, pp. 1835-1854, September 2016.
\bibitem{learning}
A.~Nedi\'{c}, A.~Olshvesky, et. al. ``Fast Convergence Rates for Distributed Non-Bayesian Learning," \emph{IEEE Trans. Autom. Control,} vol. 62, no. 11, pp. 5538-553, March 2017.
\bibitem{formation}
R.L.~Raffard, C.J.~Tomlin, et. al. ``Distributed Optimization for Cooperative Agents: Application to Formation Flight," \emph{The 43rd IEEE Annual Conference on Decision and Control,} 2004, pp. 2453-2459.
\bibitem{admm}
S.~Boyd, N.~Parikh, et. al., ``Distributed Optimization and Statistical Learning via the Alternating Direction Method of Multipliers," \emph{Foundation and Trends\textregistered in Machine Learning,} vol. 3, no. 1, pp.1-122, January 2011.
\bibitem{harnesing}
G.~Qu, N.~Li, ``Harnessing Smoothness to Accelerate Distributed Optimization," \emph{IEEE Trans. Autom. Control,} to appear.
\bibitem{application}
I.~Necoara, V.~Nedelcu, et. al., ``Parallel and distributed optimization methods for estimation and control in networks," \emph{Journal of Process Control,} vol. 21, no. 5, pp. 756-766, June 2011. 
\bibitem{nedic}
A.~Nedi\'{c}, A~Ozdaglar, ``Distributed Subgradient Methods for Multi-Agent Optimization," \emph{IEEE Trans. Autom. Control,} vol. 54, no. 1, pp. 48-60, January 2009.
\bibitem{giannakis}
G.B.~Giannakis, V.~Kekatos, et. al. ``Monitoring and Optimization for Power Grids: A Signal Processing Perspective," \emph{IEEE Signal Process. Magazine,} vol. 30, no. 5, pp. 107-128, August 2013.
\bibitem{estimation}
I.D.~Schizas, A.~Ribeiro, et. al. ``Consensus in Ad Hoc WSNs with Noisy Links-Part 1: Distributed Estimation of Deterministic Signals," \emph{IEEE Trans. Signal Process.,} vol. 56, no. 1, pp. 350-364, December 2007.
\bibitem{next}
P.~Lorenzo, G.~Scutari, ``NEXT: In-Network Nonconvex Optimization," \emph{IEEE Trans. Signal and Info. Proc. over Networks.} vol. 2, no. 2, pp. 120-136, June 2016.
\bibitem{localization}
M.G.~Rabbat, R.D.~Nowak, ``Decentralized source localization and tracking [wireless sensor networks]," \emph{Proc. IEEE ICASSP,} May 2004.
\bibitem{tutorial}
D.P.~Palomar, M~Chaing, ``A tutorial on decomposition methods for network utility maximization," \emph{IEEE Journal on Selected Areas in Communications,} vol. 8, no. 1, pp. 1439-1451, August 2006.
\bibitem{extra}
W.~Shi, Q.~Ling, et. al., ``EXTRA: An Exact First-Order Algorithm for Decentralized Consensus Optimization," \emph{SIAM J. Optim.,} vol. 25, no. 5, pp. 944-966, May 2015.
\bibitem{interpret}
A.~Nedi\'{c}, A.~Olshevsky, et. al., ``Achieving Geometric Convergence for Distributed Optimization Over Time-Varying Graphs," \texttt{arxiv:1607.03218v3}, March 2017.
\bibitem{pextra}
W.~Shi, Q.~Ling, et. al. ``A Proximal Gradient Algorithm for Decentralized Composite Optimization," \emph{IEEE Trans. Signal Process.,} vol. 63, no. 22, pp. 6013-6023, November 2015.
\bibitem{optimal}
K.~Scaman, F.~Back, et. al. ``Optimal algorithms for smooth and strongly convex distributed optimization in networks," \texttt{arxiv:1702.08704v2,} April 2017.
\bibitem{random_networks}
M.~Hong, T.~Chang, ``Stochastic Proximal Gradient Consensus over Random Networks," \emph{IEEE Trans. Signal Process.,} vol. 65, no. 11, pp. 2933-2948, June 2017.
\bibitem{time_varying_graphs}
A.~Nedi\'{c}, A.~Olshevsky, ``Distributed optimization over time-varying directed graphs," \emph{IEEE Trans. Automat. Control,} vol 60, no. 3, pp. 601-615, October 2014.
\bibitem{subgradient}
A.~Nedi\'{c} and A.~Ozdaglar, ``Distributed Subgradient Method for Multi-agent Optimization," \emph{IEEE Trans. Automat. Control,} vol 54, no. 1, pp. 48-61, January 2009.
\bibitem{dual_decomposition}
H.~Terelius, U.~Topcu, et. al, ``Decentralized multi-agent optimization via dual decomposition," \emph{Proc. IFAC,} vol. 44, no. 1, pp. 11245-11251, January 2011.
\bibitem{non_convex}
Y.~Sun, G.~Scutari, et. al. ``Distributed nonconvex Multiagent optimization Over Time-Varying networks," \emph{Asilomar}, March 2017.
\bibitem{oracle}
O.~Devolder, F.~Glineur, et. al. ``First-order methods with inexact oracle: the strongly convex case," \emph{Tech. report, UCL}, 2013.
\bibitem{dual_smooth}
S.~Kakade, S.~Shalev-Shwartz, et. al. ``On the duality of strong convexity and strong smoothness: Learning applications and matrix regularization," 2009.
\bibitem{nesterov}
Y.~Nesterov, ``Introductory lectures on Convex Optimization: A Basic Course," Springer-Verlag, US, 2004.
\bibitem{luo}
Z.Q.~Luo, P.~Tseng, ``On the Convergence Rate of Dual Ascent Methods for Linearly Constrained Convex Minimization," \emph{Math. Oper. Res.,} vol. 18, no. 4, pp. 846-867, 1993.
\bibitem{martinez}
M.~Zhu, S.~Martinez, ``Discrete-Time Dynamic Average Consensus," \emph{Automatica,} vol. 46, no. 2, pp. 322-329, 2010.
\end{thebibliography}
